\documentclass[leqno,a4paper,12pt]{article}
\usepackage{amsmath,amsthm,bbm}
\usepackage[sans]{dsfont}

\newtheorem{Thm}{\indent Theorem}[section]
\newtheorem{Prop}[Thm]{\indent Proposition}

\newtheorem{Cor}[Thm]{\indent Corollary}

\theoremstyle{definition}
\newtheorem{Def}[Thm]{\indent Definition}
\newtheorem{Rem}[Thm]{\indent Remark}
\newtheorem{Ex}[Thm]{\indent Example}
\newtheorem{Exo}[Thm]{\indent Problem}
\def\qed{{\hskip0pt\unskip\unskip\nobreak\hfil\penalty50
          \hskip1em\hbox{}\nobreak\hfil
          {\bf q.e.d.}%
          \parfillskip=0pt\finalhyphendemerits=0
          \par}\medskip}

\newenvironment{Proof}
               {{\it Proof.}\quad}
               {\qed}

\newenvironment{Proofof}[1]
               {{\it Proof of #1.}\quad}
               {\qed}

\newcommand{\Prime}{\kern3\fontdimen1\font$'$\kern-7\fontdimen1\font}

\long\def\forget#1{}

\long\def\beginSIDEREMARK#1\endSIDEREMARK
    {{\par\bigskip\advance\leftskip by 2cm
                  \advance\rightskip by -2cm\noindent
      {\bf Our own side remark:} #1
      \par\bigskip\noindent}}
\long\def\beginFORGET#1\endFORGET{#1}
\long\def\beginFORGET#1\endFORGET{}

\def\?{\ ???\ \immediate\write16{}%
\immediate\write16{Warning: There was still a question mark . . . }%
\immediate\write16{}}

\usepackage{amsmath}
\usepackage{exscale}
\usepackage{amssymb}

\input cyracc.def
\font\tencyr=wncyr6
\def\cyr{\tencyr\cyracc}
\newcommand{\cyrb}{{\cyr B}}

\newcommand{\BC}{{\mathbb{C}}}
\newcommand{\BD}{{\mathbb{D}}}

\newcommand{\BQ}{{\mathbb{Q}}}

\newcommand{\BZ}{{\mathbb{Z}}}

\newcommand{\FM}{{\mathfrak{M}}}

\newcommand{\CC}{{\cal C}}
\newcommand{\CD}{{\cal D}}

\newcommand{\CF}{{\cal F}}

\newcommand{\Spec}{\mathop{{\bf Spec}}\nolimits}

\newcommand{\imm}{\mathop{{\rm im}}\nolimits}

\newcommand{\End}{\mathop{\rm End}\nolimits}

\newcommand{\Gr}{\mathop{\rm Gr}\nolimits}

\newcommand{\Hom}{\mathop{\rm Hom}\nolimits}

\newcommand{\loccit}{[loc.$\;$cit.]}

\def\tei{\, | \,}
\def\halb{\frac{1}{2}}

\def\id{{\rm id}}

\newbox\mybox
\def\arrover#1{\mathrel{
       \setbox\mybox=\hbox spread 1.4em{\hfil$\scriptstyle#1$\hfil}
       \vbox{\offinterlineskip\copy\mybox
             \hbox to\wd\mybox{\rightarrowfill}}}}
\def\larrover#1{\mathrel{
       \setbox\mybox=\hbox spread 1.4em{\hfil$\scriptstyle#1$\hfil}
       \vbox{\offinterlineskip\copy\mybox
             \hbox to\wd\mybox{\leftarrowfill}}}}

\def\ontoover#1{\mathrel{
       \setbox\mybox=\hbox spread 1.4em{\hfil$\scriptstyle#1$\hfil}
       \vbox{\offinterlineskip\copy\mybox
             \hbox to\wd\mybox{\rightarrowfill\hskip-2.8mm
                               $\rightarrow$}}}}
\def\leftontoover#1{\mathrel{
       \setbox\mybox=\hbox spread 1.4em{\hfil$\scriptstyle#1$\hfil}
       \vbox{\offinterlineskip\copy\mybox
             \hbox to\wd\mybox{$\leftarrow$\hskip-2.8mm
                               \leftarrowfill}}}}
\def\longto{\longrightarrow}
\def\into{\hookrightarrow}

\def\longonto{\ontoover{\ }}
\def\isoto{\arrover{\sim}}

\def\longinto{\lhook\joinrel\longrightarrow}

\usepackage[curve,matrix,arrow,cmtip]{xy}
\NoComputerModernTips

\def\myxymessage{\def\messagetext
   {Here an xy-pic diagram was omitted to speed up compilation . . . }
   \immediate\write16{\messagetext}
   \hbox{\bf \messagetext}}
\def\filxymatrix#1{\myxymessage}
\def\filxyarray#1{\myxymessage}
\newdir^{ (}{{}*!/-3pt/\dir^{(}}
\newdir_{ (}{{}*!/-3pt/\dir_{(}}
\newdir^{ )}{{}*!/+3pt/\dir^{)}}
\newdir_{ )}{{}*!/+3pt/\dir_{)}}

\def\rscript#1{\hbox to 0pt{$\scriptstyle#1$\hss}}

\let\oldbullet\bullet
\def\bullet{{\mathchoice{\oldbullet}%
                        {\oldbullet}%
                        {\scriptscriptstyle\oldbullet}%
                        {\oldbullet}}}

\newcommand{\argdot}{{\;\bullet\;}}%Punkt als Platzhalter fuer Argumente
\newcommand{\nX}{\mathop{X^{norm}}\nolimits}
\newcommand{\nnX}{\mathop{X^{nn}}\nolimits}
\newcommand{\Xp}{\mathop{\widetilde{X}}\nolimits}
\newcommand{\nj}{\mathop{j^{norm}}\nolimits}
\newcommand{\ip}{\tilde{\imath}}
\newcommand{\nni}{\mathop{i^{nn}}\nolimits}

\newcommand{\CHUM}{\mathop{CHM(U)_\BQ}\nolimits}
\newcommand{\CHXM}{\mathop{CHM(X)_\BQ}\nolimits}
\newcommand{\CHZM}{\mathop{CHM(Z)_\BQ}\nolimits}

\newcommand{\DBM}{\mathop{DM_{\text{\cyrb}}}\nolimits}
\newcommand{\DBcM}{\mathop{DM_{\text{\cyrb},c}}\nolimits}

\newcommand{\DBXM}{\mathop{\DBM(X)}\nolimits}
\newcommand{\DBXrM}{\mathop{\DBM(X_{red})}\nolimits}
\newcommand{\DBYM}{\mathop{\DBM(Y)}\nolimits}
\newcommand{\DBZM}{\mathop{\DBM(Z)}\nolimits}
\newcommand{\DBcXM}{\mathop{\DBcM(X)}\nolimits}
\newcommand{\MHS}{\mathop{\bf MHS}\nolimits}
\newcommand{\MHM}{\mathop{\bf MHM}\nolimits}

\newcommand{\one}{\mathds{1}}

\newcommand{\uHom}{\mathop{\underline{Hom}}\nolimits}

\begin{document}

%%%%%%%%%%%%%%%%%%%%%%%%%%%%%%%%%%%%%%%%%%%%%%%%%%%%%%%%%%%%%%%%%%%%%%%
%
%  formatting

\hfuzz=3pt
\overfullrule=10pt                   % erzeugt schwarze Fehlerbalken

% The displayskip values were changed because \LaTeX does not react
% correctly to a \leqno: it should then use big skips, but doesn't.

\setlength{\abovedisplayskip}{6.0pt plus 3.0pt}
                               % preset 10.0pt plus 2.0pt minus 5.0pt
\setlength{\belowdisplayskip}{6.0pt plus 3.0pt}
                               % preset 10.0pt plus 2.0pt minus 5.0pt
\setlength{\abovedisplayshortskip}{6.0pt plus 3.0pt}
                               % preset 0.0pt plus 3.0pt
\setlength{\belowdisplayshortskip}{6.0pt plus 3.0pt}
                               % preset 6.0pt plus 3.0pt minus 3.0pt

\setlength{\baselineskip}{13.0pt}
                               % preset 12.0pt
\setlength{\lineskip}{0.0pt}
                               % preset 1.0pt
\setlength{\lineskiplimit}{0.0pt}
                               % preset 0.0pt

%%%%%%%%%%%%%%%%%%%%%%%%%%%%%%%%%%%%%%%%%%%%%%%%%%%%%%%%%%%%%%%%%%%%%%%
%
%  Title Page
%
%%%%%%%%%%%%%%%%%%%%%%%%%%%%%%%%%%%%%%%%%%%%%%%%%%%%%%%%%%%%%%%%%%%%%%%

\title{Motivic intersection complex
\forget{
\footnotemark
\footnotetext{To appear in ....}
}
}
\author{\footnotesize by\\ \\
\mbox{\hskip-2cm
\begin{minipage}{6cm} \begin{center} \begin{tabular}{c}
J\"org Wildeshaus \footnote{
Partially supported by the \emph{Agence Nationale de la
Recherche}, project no.\ ANR-07-BLAN-0142 ``M\'ethodes \`a la
Voevodsky, motifs mixtes et G\'eom\'etrie d'Arakelov''. }\\[0.2cm]
\footnotesize LAGA\\[-3pt]
\footnotesize UMR~7539\\[-3pt]
\footnotesize Institut Galil\'ee\\[-3pt]
\footnotesize Universit\'e Paris 13\\[-3pt]
\footnotesize Avenue Jean-Baptiste Cl\'ement\\[-3pt]
\footnotesize F-93430 Villetaneuse\\[-3pt]
\footnotesize France\\
{\footnotesize \tt wildesh@math.univ-paris13.fr}
\end{tabular} \end{center} \end{minipage}
\hskip-2cm}
\\[2.5cm]
\forget{
{\bf Preliminary version --- not for distribution}\\[1cm]
}
}
% In the final version we might want to fix the date:
\date{April 28, 2011}
\maketitle
%\quad \\[-1.7cm]
\begin{abstract}
\noindent 
In this article, we give an unconditional
definition of the motivic analogue of the intersection complex, 
establish its basic properties, and prove its existence in certain cases.\\


%\noindent
%{\bf R\'esum\'e~:} \\
\end{abstract} 
%\vfill

\bigskip
\bigskip
\bigskip

\noindent {\footnotesize Math.\ Subj.\ Class.\ (2010) numbers: 19E15
(14C25, 14F42, 14J17).
}

\eject
\tableofcontents

\bigskip
\vspace*{0.5cm}

%\newpage
%\include{Intro}
%%%%%%%%%%%%%%%%%%%%%%%%%%%%%%%%%%%%%%%%%%%%%%%%%%%%%%%%%%%%%%%%%%%%%%%
%
%  Introduction
%
%%%%%%%%%%%%%%%%%%%%%%%%%%%%%%%%%%%%%%%%%%%%%%%%%%%%%%%%%%%%%%%%%%%%%%%

\setcounter{section}{-1}
\section{Introduction}
\label{Intro}

%%%%%%%%%%%%%%%%%%%%%%%%%%%%%%%%

%%%%%%%%%%%%%%%%%%%%%%%%%%%%%%%%

This paper contains largely extended notes of the talk the author gave
during the conference \emph{Regulators~III}, which took place at
the University of Barcelona in July 2010. 
Its main purpose is to propose an unconditional definition
of the motivic intersection complex. \\

The use of the intersection complex, say in the context of (topological) sheaves
on schemes over the complex numbers, or of ($\ell$-adic) sheaves on schemes
over a field, can be motivated by \emph{purity}. Let $X$ be proper over $k$. 
Its singular cohomology (if $k = \BC$) carries
a pure Hodge structure, and its $\ell$-adic cohomology (if $k$ is finite
or a number field) a pure Galois action, provided that $X$ is smooth.
If this latter hypothesis is not met, then in order to get analogous purity statements, 
the constant sheaf on $X$ has to be replaced by the intersection complex \cite{BBD} 
(with respect to the inclusion of the regular locus of $X$). 
Its (hyper)cohomology is known as \emph{intersection cohomology} of $X$. \\

One of the main arithmetic applications to keep in mind concerns the 
\emph{Baily--Borel compactification}
of a smooth \emph{Shimura variety}: it is canonical, and even minimal in a precise sense,
but rarely smooth. Its intersection cohomology contains valuable arithmetic
information, \emph{e.g.}, certain of its direct factors allow to realize Hodge  
structures and Galois representations associated to automorphic forms. \\

In order to construct motives inducing these Hodge structures and Galois representations
\emph{via} the respective realizations, one is thus led to try first to construct
the \emph{intersection motive}. One minimal requirement on this object would be
that its realizations equal intersection cohomology. \\ 

This construction succeeded in a small number of cases. Let us cite 
varieties (over $\BC$) admitting \emph{semismall resolutions} \cite{CM}, which includes
the case of surfaces, and Baily--Borel compactifications
of Hilbert--Blumenthal varieties \cite{GHM} (over $\BC$, and with more general 
than just constant coefficients). 
A general program for the construction of the intersection motive, assuming
Grothendieck's standard conjectures, was developed (still over $\BC$) in \cite{CH}. \\ 

When the construction works unconditionally, then it does so for specific geometric reasons.
For example, such a reason would be that the relevant cycle classes are isomorphisms.
The idea is basically to obtain an explicit formula for intersection cohomology
sitting in the cohomology of a desingularization of $X$;
the specific geometric reasons in question then allow to give a motivic sense to
the explicit formula. Unfortunately, some functoriality properties
valid for intersection cohomology are not a consequence of the explicit
formula, and hence do not obviously hold for the intersection motive.
This concerns for example the action of the Hecke algebra (which is
needed in order to cut out the motive of an individual automorphic form
from the intersection motive). \\ 

In \cite{W1}, we gave an unconditional construction of the intersection motive
of Baily--Borel compactifications of smooth Hilbert--Blumenthal varieties
with non-constant coefficients. It is a Chow motive over $\BQ$,
and behaves well under Hecke correspondences. 
Again, the construction works for specific geometric reasons, which translate into saying
that ``the boundary avoids weights $-1$ and $0$''. Let us not worry
about the precise meaning of the ``boundary'' here. Rather, let us concentrate
on the central notion of \emph{weight}. \\

Assume first that our base scheme $X$ equals the spectrum of a perfect field $k$. 
According to Bondarko \cite{Bo1}, the category of \emph{geometrical motives}
\cite{VSF} carries a \emph{weight structure}, whose \emph{heart}
equals the category $CHM(k)$ of Chow motives over $k$. The precise definitions
of weight structures and hearts will be recalled in the present Section~\ref{1}; 
for the moment, let us keep in mind that the motivic weight structure allows 
for an \emph{intrinsic} characterization of the full sub-category
$CHM(k)$ of the category of geometrical motives. This is the key
for everything to follow. Roughly speaking,
the construction from \cite{W1} works since geometrical motives
are flexible enough to preserve functoriality; the problem of knowing
whether the result of this construction 
is a Chow motive is then reduced to a computation of weights. \\

In general, the properties of intersection cohomology (functoriality,
purity,...) are consequences of properties of the intersection complex.
A general solution to the problem of constructing the intersection motive
therefore requires the construction of the motivic
intersection complex. 
Here, one is confronted with a foundational problem: the na\"{\i}ve gene\-ralization
of the de\-fi\-nition via truncations \cite{BBD} cannot work since it requires the existence 
of a (perverse)
\emph{$t$-structure}. But even when the base is of the form $\Spec k$, then 
except for certain fields $k$, such a $t$-structure is not known to exist 
on the category of geometrical motives. Thus, the mere problem of giving an unconditional
definition of the motivic intersection complex is \emph{a priori} non-trivial. \\

The solution to this problem that we shall propose, is again based on the notion
of weight structure. In a way, our approach can be seen as 
``reading \cite{BBD} backwards'', \emph{i.e.}, starting from \cite[Chap.~5]{BBD} 
on weights. This concerns in particular 
the Decomposition Theorem \cite[Thm.~5.4.5]{BBD}, which implies that every
pure complex on $X$ restricting to the structure sheaf on an open smooth
sub-variety, contains the intersection complex as a direct factor.
Let us indicate already here that the motivic analogue of this result 
(Theorem~\ref{3A}~(b)) is a rather elementary exercice in weight structures... \\

Let us now give a detailed overview of the individual sections of this paper.
Section~\ref{1} starts with a review of \emph{Beilinson motives} \cite{CD}, which
conveniently generalize geometrical motives from $\Spec k$
to arbitrary bases $X$. We then recall the basic notions related to weight structures.
We review the main results from \cite{H} on the existence of the motivic weight
structure on Beilinson motives (generalizing \cite{Bo1} from $\Spec k$ to $X$), and on the behaviour
of weights under the six operations from \cite{CD}. We then \emph{define}
the category $\CHXM$ of \emph{Chow motives} over $X$ as the heart of the motivic weight structure,
and establish two complements of the theory. First (Theorem~\ref{1F}), we show that
for an open sub-scheme $U$ of $X$, the inverse image from $\CHXM$ to $\CHUM$
is both essentially surjective and full. Following the terminology introduced in \cite{Bo2}, 
this can be seen as a \emph{motivic version of resolution of singularities}.
Theorem~\ref{1F} strenghens \cite[Thm.~2.2.1~III~1]{Bo2}, where essential surjectivity is proved
up to pseudo-Abelian completion. Second (Theorem~\ref{1G}), we show that local duality
respects the weights in a strict sense; in particular, the dual of a Chow motive
is again a Chow motive. This complements \cite[Cor.~2.2.5]{H2}, where the same result is proved
provided $X$ is regular, and also \cite[Cor.~3.9]{H},
where left exactness (with respect to the weights) is established 
for any $X$. \\

Having in mind the Decomposition Theorem
\cite[Thm.~5.4.5]{BBD}, an \emph{intermediate extension} of a Chow motive 
$M_U$ over a dense $U$ should satisfy a certain minimality condition
among all possible extensions of $M_U$ to 
a Chow motive over $X$. In Section~\ref{2}, we make this precise for regular $U$,
and $M_U = \one_U$, the structure motive on $U$. 
More precisely (Definition~\ref{2A}), the motivic intersection
complex $j_{!*} \one_U$ is a Chow motive on $X$ restricting to give $\one_U$, and admitting no 
non-trivial endomorphisms restricting trivially to $U$. We then
establish independence of $j_{!*} \one_U$ of $U$ (Proposition~\ref{2B}). 
In its essence, it results from the study of a basic, but important example: 
when $X$ is regular, then $j_{!*} \one_U = \one_X$ (Example~\ref{2ex}). \\ 

Section~\ref{3} contains our main results.
According to Theorem~\ref{3A}~(a), the motivic intersection complex is unique
up to unique isomorphism. As already indicated, Theorem~\ref{3A}~(b) states that
any extension of $\one_U$ to a Chow motive over $X$ contains $j_{!*} \one_U$
as a direct factor --- provided the latter exists. Under the same hypothesis,
$j_{!*} \one_U$ is auto-dual (Corollary~\ref{3C}), meaning that the motivic
intersection pairing can be defined. Theorem~\ref{3Z} identifies the few
cases where we actually know the motivic intersection complex to exist. 
Section~\ref{4} contains the proof of Theorem~\ref{3Z}. \\

We choose to add a number of ``Problems'' in the text.
While they concern properties that one might reasonably 
expect $j_{!*} \one_U$ to satisfy, the author 
does not know to solve any of them. 
The paper also contains a number of miscellaneous
results, which are not needed elsewhere in the text, but seem worth to be mentioned nonetheless.
In particular, this concerns Corollaries~\ref{1Bei} and \ref{3loc}.
The first (Corollary~\ref{1Bei}) states that for an open immersion $j: U \into X$ 
and any Chow motive $N_U$
over $U$, the image under the inverse image $j^*$ of motivic cohomology of $N$
in motivic cohomology of $N_U$ is independent of the extension of $N_U$ to a Chow motive
$N$ over $X$. We relate this to Scholl's construction of ``integral'' sub-spaces of
motivic cohomology of Chow motives over number fields (Remark~\ref{1Sch}).
According to the second (Corollary~\ref{3loc}), a Beilinson motive which is 
Nisnevich-locally isomorphic to $\one_X$, for a regular base $X$, is (globally)
isomorphic to $\one_X$. This allows to generalize \emph{absolute purity}
\cite[Thm.~13.4.1]{CD} to arbitrary morphisms $a: X \to S'$ between regular schemes: as soon
as $a$ is of pure relative dimension $d$, there is an isomorphism
$\one_X (d)[2d] \cong a^! \one_{S'}$  
(Corollary~\ref{3abs}). \\

Part of this work was done while I was enjoying a 
\emph{modulation de service pour les porteurs de projets de recherche},
granted by the \emph{Universit{\'e} Paris~13},
and during a stay at the University of Tokyo.
I am grateful to both institutions.
I wish to thank the organizers of \emph{Regulators~III}
for the invitation to Barcelona, C.~Soul\'e for a stimulating 
question asked during my talk, F.~D\'eglise and D.~H\'ebert for useful comments
on a first draft of this paper, G.~Ancona for strengthening 
an earlier version of Proposition~\ref{2C}, and the referee for 
her or his remarks and suggestions. \\

{\bf Conventions}: Throughout the article, 
$S$ denotes a fixed base scheme, which we assume
to be of finite type over an excellent scheme 
of dimension at most two. By definition, \emph{schemes} are  
$S$-schemes which are separated and 
of finite type (in particular, they are excellent, and
Noetherian of finite dimension), \emph{morphisms} between schemes
are separated morphisms of $S$-schemes, and 
a scheme is \emph{regular}
if the underlying reduced scheme is regular in the usual sense.

%%% Local Variables:
%%% mode: latex
%%% TeX-master: "head"
%%% End:

\bigskip

%\include{Sec1}
%%%%%%%%%%%%%%%%%%%%%%%%%%%%%%%%%%%%%%%%%%%%%%%%%%%%%%%%%%%%%%%%%%%%%%%
%
%  Section 1
%
%%%%%%%%%%%%%%%%%%%%%%%%%%%%%%%%%%%%%%%%%%%%%%%%%%%%%%%%%%%%%%%%%%%%%%%

\section{Review of weights on Beilinson motives}
\label{1}

%%%%%%%%%%%%%%%%%%%%%%%%%%%%%%%%

%%%%%%%%%%%%%%%%%%%%%%%%%%%%%%%%

We fix our base $S$, and work in the triangulated, $\BQ$-linear categories
$\DBXM$ of \emph{Beilinson motives} over $X$ \cite[Def.~13.2.1]{CD}, 
indexed by schemes $X$ (always in the sense of the conventions
fixed at the end of our Introduction). 
As in \cite{CD}, the symbol $\one_X$
is used to denote the unit for the tensor product in $\DBXM$.
We shall employ the full formalism of six operations developed in
\loccit . Below, we shall 
list the principles (A)--(E) which will be parti\-cu\-lar\-ly important
to us. The global assumptions made in \loccit \
to establish these principles are met since $\DBM (\argdot)$
is a \emph{motivic category} \cite[Cor.~13.2.11]{CD}, which by 
definition \cite[Def.~2.4.2]{CD} implies that it is
\emph{pregeome\-tric}. Therefore, it is \emph{$Sm$-fibred}
\cite[Def.~1.1.9]{CD},
the \emph{localization property} $(Loc_i)$ from \cite[Def.~2.3.2]{CD} holds,
and by \cite[Thm.~2.4.12]{CD} so do the 
\emph{proper transversality property} from \cite[Def.~1.1.16]{CD} and the
\emph{support property} from \cite[Def.~2.2.5]{CD}.
Furthermore, by \cite[Prop.~14.2.16]{CD},
the category $\DBM (\argdot)$ is \emph{separated} 
in the sense of \cite[Def.~2.1.11]{CD}. By \cite[Ex.~14.3.20]{CD},
it is \emph{pure} in the sense of \cite[Def.~14.3.19]{CD}.
(A)~\emph{Absolute purity. Relation to $K$-theory:} if $i: Z \into X$ 
is a closed immersion of pure codimension $c$ between regular schemes,
then there is a canonical isomorphism
\[  
\one_Z(-c)[-2c] \isoto i^! \one_X 
\]
in $\DBZM$ \cite[Thm.~13.4.1]{CD}. 
For any regular scheme $X$, and
any pair of integers $(p,q)$, there is a canonical isomorphism
\[
\Hom_{\DBXM}(\one_X , \one_X(p)[q]) \cong \Gr_\gamma^p K_{2p-q}(X)_\BQ \; ,
\]
where $K_\bullet(X)_\BQ$ denotes the tensor product of $K$-theory of $X$
with the rationals, and $\Gr_\gamma$ the graded object with respect
to the (Adams) gamma filtration \cite[Cor.~13.2.14]{CD}. 
Furthermore, this isomorphism is contravariantly functorial with respect
to morphisms of regular schemes \cite[Cor.~13.2.11]{CD}.
(B)~\emph{Base change:} for any morphism $f$,
there is a natural transformation
\[
\alpha_f : f_! \longto f_* \; ,
\]
which is an isomorphism is $f$ is proper \cite[Thm.~2.2.14~(1)]{CD}. 
If $f$ is the base of a cartesian diagram
\[
\vcenter{\xymatrix@R-10pt{
        Y' \ar[d]_{g'} \ar[r]^{f'}  & 
        X' \ar[d]^g \\
        Y \ar[r]^f &
        X  
\\}}
\]
of schemes, then the exchange transformation
\[
g^*f_! \longto f'_!{g'}^*
\]
is an isomorphism \cite[Prop.~2.2.13~(b)]{CD}. Hence so is the 
adjoint exchange transformation
\[
g'_*{f'}^! \longto f^!g_* \; .
\] 
(C)~\emph{Constructibility:} by definition \cite[Def.~1.4.7]{CD},
the full thick triangulated sub-category $\DBcXM$ of $\DBXM$
of constructible objects is generated by the Tate twists $M_X(T)(p)$
of the \emph{motives} $M_X(T)$ \cite[Sect.~1.1.33]{CD} of smooth 
$X$-schemes $T$. In particular, all twists $\one_X(p)$ belong
to $\DBcXM$. By \cite[Ex.~14.1.3]{CD}, an object of $\DBXM$
is constructible if and only if it is compact. 
According to 
\cite[Thm.~14.1.31]{CD}, the sub-categories
$\DBcM(\argdot) \subset \DBM(\argdot)$ are respected by the six functors.
(D)~\emph{Duality:} fix a scheme $X$ whose structure morphism 
to $S$ factors over a regular scheme; this is of course the case 
if $S$ is itself regular.
According to \cite[Thm.~14.3.28]{CD},
the category $\DBcXM$ then contains \emph{dualizing objects} 
in the sense of \cite[Def.~14.3.10]{CD}. 
Fix such a dualizing object $R$.
Define the \emph{local duality functor} (with respect to $R$) as 
\[
\BD_X := \uHom_X (\argdot , R) \; .
\]
It is right adjoint to itself \cite[Sect.~14.3.30]{CD}.
It preserves constructible objects, and the adjunction
$\id_X \to \BD_X^2$ is an isomorphism on $\DBcXM$
\cite[Cor.~14.3.31~(a), (b)]{CD}. Furthermore, it exchanges 
$f^*$ and $f^!$, as well as $f_!$ and $f_*$ in the following sense:
for a morphism $f:Y \to X$, put 
\[
\BD_Y := \uHom_Y (\argdot , f^! R) \; ;
\]  
note that according to \cite[Prop.~14.3.29~(ii)]{CD}, the motive
$f^! R$ is dualizing on $Y$.
Then there are natural isomorpisms of functors 
\[
\BD_Y f^* \isoto f^! \BD_X \quad \text{and} \quad
f_* \BD_Y \isoto \BD_X f_! 
\]
on $\DBM(\argdot)$
\cite[Cor.~14.3.31~(d) and its proof]{CD}. Therefore,
\[
f^* \BD_X \isoto f^! \BD_Y \quad \text{and} \quad
\BD_X f_* \isoto f_! \BD_Y 
\]
on $\DBcM(\argdot)$.
For the applications of duality that we have in mind,
we need to make explicit choices of 
dualizing object $R$. Fix a pair of integers $(p,q)$,
and a morphism $a: X \to S'$ with regular target. Then
\[
R := a^! \one_{S'}(p)[q] \in \DBcXM 
\]
is a dualizing object \cite[Prop.~14.3.29]{CD}.
It will be necessary to identify $R$ under the following additional
hypotheses on $X$: the morphism $a: X \to S'$ is quasi-projective,
and $X$ is regular and connected of relative dimension $e$ over $S'$.
We claim that in this case, there is an isomorphism
\[
R \isoto \one_X(p+e)[q+2e] \; .
\]
Indeed, absolute purity (see point~(A)) and the formula $j^! = j^*$
for an open immersion $j$ \cite[Thm.~2.2.14~(2)]{CD} reduces us
to the case when $X$ is a projective space over $S'$. 
Our claim then follows from
\cite[Scholie~1.4.2~3]{A} (\emph{via} \cite[Cor.~2.4.9]{CD}).
(E)~\emph{Localization:} if $i: Z \into X$
and $j: U \into X$ are complementary closed, resp.\ open immersions
of schemes, then there are canonical exact triangles 
\[
j_!j^* \longto \id_X \longto i_*i^* \longto j_!j^*[1] \; ,
\]
\[
i_*i^! \longto \id_X \longto j_*j^* \longto i_*i^![1] 
\]
of exact endo-functors of $\DBXM$ 
\cite[Prop.~2.3.3~(2), (3), Thm.~2.2.14~(2)]{CD}. 
The adjunctions $\id_U \to j^*j_!$, $j^*j_* \to \id_U$ and $i^*i_* \to \id_Z$ 
are isomorphisms, and the compositions $i^* j_!$
and $j^* i_*$ are trivial \cite[Sect.~2.3.1]{CD}. 
From what precedes, it follows formally that the adjunction $\id_Z \to i^!i_*$
is an isomorphism, and that the composition $i^!j_*$ is trivial. 
We also see, putting $i$ equal to the
immersion of the 
reduced scheme structure $X_{red}$ on $X$, that 
\[
i_* : \DBXrM \longto \DBXM
\]
is an equivalence of categories, with canonical quasi-inverse $i^! = i^*$. 
This jusitifies \emph{a posteriori} the abuse of language fixed 
in the conventions at the end of our Introduction. \\

Now recall the following notions, due to Bondarko.

\begin{Def}[{\cite[Def.~1.1.1]{Bo1}}] \label{1A}
Let $\CC$ be a triangulated category. A \emph{weight structure on $\CC$}
is a pair $w = (\CC_{w \le 0} , \CC_{w \ge 0})$ of full 
sub-categories of $\CC$, such that, putting
\[
\CC_{w \le n} := \CC_{w \le 0}[n] \quad , \quad
\CC_{w \ge n} := \CC_{w \ge 0}[n] \quad \forall \; n \in \BZ \; ,
\]
the following conditions are satisfied.
\begin{enumerate}
\item[(1)] The categories
$\CC_{w \le 0}$ and $\CC_{w \ge 0}$ are 
Karoubi-closed: for any object $M$ of $\CC_{w \le 0}$ or
$\CC_{w \ge 0}$, any direct summand of $M$ formed in $\CC$
is an object of $\CC_{w \le 0}$ or
$\CC_{w \ge 0}$, respectively.
\item[(2)] (Semi-invariance with respect to shifts.)
We have the inclusions
\[
\CC_{w \le 0} \subset \CC_{w \le 1} \quad , \quad
\CC_{w \ge 0} \supset \CC_{w \ge 1}
\]
of full sub-categories of $\CC$.
\item[(3)] (Orthogonality.)
For any pair of objects $A \in \CC_{w \le 0}$ and $B \in \CC_{w \ge 1}$,
we have
\[
\Hom_{\CC}(A,B) = 0 \; .
\]
\item[(4)] (Weight filtration.)
For any object $M \in \CC$, there exists an exact triangle
\[
A \longto M \longto B \longto A[1]
\]
in $\CC$, such that $A \in \CC_{w \le 0}$ and $B \in \CC_{w \ge 1}$.
\end{enumerate}
\end{Def}

Slightly generalizing the above terminology, for $n \in \BZ$,
we shall refer to any exact triangle
\[
A \longto M \longto B \longto A[1]
\]
in $\CC$, with $A \in \CC_{w \le n}$ and $B \in \CC_{w \ge n+1}$,
as a weight filtration of $M$.

\begin{Def}[{\cite[Def.~1.2.1~1]{Bo1}}] \label{1B}
Let $w$ be a weight structure on $\CC$.
The \emph{heart of $w$} is the full additive sub-category $\CC_{w = 0}$
of $\CC$ whose objects lie 
both in $\CC_{w \le 0}$ and in $\CC_{w \ge 0}$.
\end{Def}

Beilinson motives can be endowed with weight structures,
thanks to the main results from \cite{H}. More precisely, the following holds.

\begin{Thm}[{\cite[Thm.~3.3, Thm.~3.8~(i)--(ii)]{H}}] \label{1C}
(a)~There are ca\-no\-ni\-cal weight structures $w$ on the categories 
$\DBcM(\argdot)$. They are uniquely characterized by the following properties. 
\begin{enumerate}
\item[(a1)] The objects $\one_X(p)[2p]$ belong
to the heart $\DBcM(X)_{w = 0}$, for all integers $p$, whenever $X$ is regular.
\item[(a2)] For a morphism of schemes $f$, left adjoint functors 
$f^*$, $f_!$ and $f_\sharp$ (the latter for smooth $f$) are 
\emph{$w$-left exact}, \emph{i.e.}, they map
$\DBcM(\argdot)_{w \le 0}$ to $\DBcM(\argdot)_{w \le 0}$, and right adjoint
functors $f_*$, $f^!$ and $f^*$  (the latter for smooth $f$) are
\emph{$w$-right exact}, \emph{i.e.}, they map
$\DBcM(\argdot)_{w \ge 0}$ to $\DBcM(\argdot)_{w \ge 0}$.
\end{enumerate}
(b)~There are canonical weight structures $W$ on the categories 
$\DBM(\argdot)$. They induce the weight structures $w$ on the categories
$\DBcM(\argdot)$. They are uniquely characterized by the requirement 
that any small sum of objects of $\DBcM(X)_{w = 0}$ 
lie in $\DBM(X)_{W = 0}$.
\end{Thm}

Let us refer to the weight structure $w$ on $\DBcM(\argdot)$ as the
\emph{motivic weight structure}.
Theorem~\ref{1C} generalizes an earlier result 
of Bondarko's \cite[Prop.~6.5.3]{Bo1} concerning the case
$X = \Spec k$, for a perfect field $k$ (use
\cite[Rem.~10.1.5, Thm.~15.1.4]{CD} to get the equivalence
between the triangulated category of \emph{geometrical motives}
\`a la Voevodsky and $\DBcM(\Spec k)$). 

\begin{Rem}
\forget{
(a)~The reader may choose to work in the categories 
denoted ${\bf SH}_\FM(\argdot)$
of \cite[D\'ef.~4.5.52]{A} instead of $\DBM(\argdot)$, with $\FM$
equal to the category of complexes of vector spaces over $\BQ$,
and $\tau$ equal to the \'etale topology \cite[beginning of Sect.~4.5]{A}.
According to \cite[Thm.~15.2.16]{CD}, the two points of view are equivalent.
In \cite{A}, the formalism of six functors had been developed intrinsically
for quasi-projective morphisms of schemes, which is sufficient for the
purposes of the present paper since all morphisms to be considered
are indeed quasi-projective. \\[0.1cm]
(b)~
}Since the first appearance of \cite{H}, a different proof of
existence of the motivic weight structure 
was given in \cite[Thm.~2.1.1]{Bo2}. The $w$-exactness
properties from \cite[Thm.~3.8]{H} are shown in \cite[Thm.~2.2.1~II]{Bo2} 
for quasi-projective morphisms of schemes. The results of \cite{Bo2} 
were obtained independently from \cite{H}. 
\end{Rem}

Note that for perfect fields $k$,
\cite[Sect.~6.6]{Bo1} allows to identify the heart
of the motivic weight structure on $\DBcM(\Spec k)$
with the category (opposite to the category)
of Chow motives over $k$. This motivates the following.

\begin{Def} \label{1D}
The $\BQ$-linear category $\CHXM$ of \emph{Chow motives} over $X$
is defined as the heart $\DBcM(X)_{w = 0}$ of the motivic weight structure.
\end{Def}

\begin{Rem} \label{1E}
The categories $\DBcM(\argdot)$ are pseudo-Abelian (see \cite[Sect.~2.10]{H}). 
Hence so are their hearts $CHM(\argdot)_\BQ$.
For a fixed scheme $X$, the category $\CHXM$
can be constructed as the pseudo-Abelian completion of the category
of motives over $X$ of the form
\[
f_! \one_Y (p)[2p] \; ,
\]
for proper morphisms $f: Y \to X$ with regular source $Y$, and integers $p$
\cite[Thm.~3.3~(ii)]{H}. Since by \cite[Cor.~14.3.9]{CD} these motives
generate $\DBcXM$ as a thick triangulated category, we see in particular that
the latter is generated by the heart of its weight structure. 
\end{Rem}

Here is our first application of the formalism of motivic weight structures.

\begin{Thm} \label{1F}
Let $j: U \into X$ be an open immersion of schemes. \\[0.1cm]
(a)~The inverse image
\[
j^*: \CHXM \longto \CHUM
\]
is essentially surjective. \\[0.1cm]
(b)~The inverse image $j^*$ is full.
\end{Thm}

Note that by Theorem~\ref{1C}~(a2), the functor $j^*$ is \emph{$w$-exact},
meaning that it is both $w$-left
and $w$-right exact ($j$ is smooth). In particular, it preserves 
the hearts of the weight structures on $\DBcM(X)$ and on $\DBcM(U)$.
Note also that essential surjectivity of $j^*$ on both 
$\DBcM(\argdot)_{w \le 0}$ and $\DBcM(\argdot)_{w \ge 0}$ is a formal
consequence of the existence of $j_!$ and $j_*$, and the formulae
$\id_U \cong j^*j_!$ and $j^*j_* \cong \id_U$.
(By contrast, $j^*$ should not in general be expected to be full
on $\DBcM(\argdot)_{w \le 0}$ or on $\DBcM(\argdot)_{w \ge 0} \, $!)
Theorem~\ref{1F}~(a) 
strenghens \cite[Thm.~2.2.1~III~1]{Bo2}, where it is proved 
that $j^*$ is essentially surjective up to pseudo-Abelian completion.

\medskip

\begin{Proofof}{Theorem~\ref{1F}}
(a)~Let $M_U$ be an object of $\CHUM$, and consider the morphism
\[
m:= \alpha_j(M_U): j_! M_U \longto j_* M_U
\]
of motives over $X$ (see point~(B) above). 
Applying $j^*$ to $m$ yields an isomorphism.
Therefore, by localization, 
any cone of $m$ is of the form $i_* C$,
for a motive $C$ over the complement $i: Z \into X$ of $U$ in $X$
(with the reduced scheme structure, say).
Choose and fix such a cone $i_* C$, as well as a weight filtration
\[
C_{\le 0} \stackrel{c_-}{\longto} C 
\stackrel{c_+}{\longto} C_{\ge 1} 
\stackrel{\delta}{\longto} C_{\le 0}[1]
\]
of $C \in \DBcM(Z)$ (Theorem~\ref{1C}~(a)). 
Thus, 
\[
C_{\le 0} \in \DBcM(Z)_{w \le 0} \quad \text{and} \quad 
C_{\ge 1} \in \DBcM(Z)_{w \ge 1} \; .
\] 
According to axiom TR4' of triangulated categories (see \cite[Sect.~1.1.6]{BBD}
for an equivalent formulation), the diagram of exact triangles
\[
\vcenter{\xymatrix@R-10pt{
        0 \ar[d] \ar[r] &
        i_* C_{\ge 1}[-1] \ar@{=}[r] &
        i_* C_{\ge 1}[-1] \ar[d]^{i_* \delta[-1]} \ar[r] &
        0 \ar[d] \\
        j_! M_U \ar@{=}[d] &  
         &
        i_* C_{\le 0} \ar[d]^{i_* c_-} \ar[r] &
        j_! M_U[1] \ar@{=}[d] \\
        j_! M_U \ar[d] \ar[r]^m &
        j_* M_U \ar[d] \ar[r] &
        i_* C \ar[d]^{i_* c_+} \ar[r] &
        j_! M_U[1] \ar[d] \\
        0 \ar[r] &
        i_* C_{\ge 1} \ar@{=}[r] &
        i_* C_{\ge 1} \ar[r] &
        0    
\\}}
\]
in $\DBcM(X)$ can be completed to give
\[
\vcenter{\xymatrix@R-10pt{
        0 \ar[d] \ar[r] &
        i_* C_{\ge 1}[-1] \ar[d] \ar@{=}[r] &
        i_* C_{\ge 1}[-1] \ar[d]^{i_* \delta[-1]} \ar[r] &
        0 \ar[d] \\
        j_! M_U \ar@{=}[d] \ar[r] &  
        M \ar[d] \ar[r] &
        i_* C_{\le 0} \ar[d]^{i_* c_-} \ar[r] &
        j_! M_U[1] \ar@{=}[d] \\
        j_! M_U \ar[d] \ar[r]^m &
        j_* M_U \ar[d] \ar[r] &
        i_* C \ar[d]^{i_* c_+} \ar[r] &
        j_! M_U[1] \ar[d] \\
        0 \ar[r] &
        i_* C_{\ge 1} \ar@{=}[r] &
        i_* C_{\ge 1} \ar[r] &
        0    
\\}}
\] 
with $M \in \DBcM(X)$. Since the composition of functors $j^* i_*$ 
is trivial, the inverse image $j^* M$ is
isomorphic to $M_U$. Now observe that by Theorem~\ref{1C}~(a2),
the functors $i_! = i_*$ 
and $j_!$ are $w$-left exact, and $i_*$ and $j_*$ are 
$w$-right exact. Thus, by the above diagram, the motive $M$
is simultaneously an extension of motives of weights $\le 0$, 
and an extension of motives of weights $\ge 0$.
It follows easily (see \cite[Prop.~1.3.3~3]{Bo1})
that $M$ is pure of weight zero.

(b)~Now let $M$ and $N$ be Chow motives over $X$, and assume that 
a morphism
\[
\beta_U : j^* M \longto j^* N
\]
between their restrictions to $U$ is given.
Consider the localization triangles for
$M$ and for $N$.
\[
\vcenter{\xymatrix@R-10pt{
        i_* i^* M [-1] \ar[r] &
        j_! j^* M \ar[d]^{j_! \beta_U} \ar[r] &
        M \ar[r] &
        i_* i^* M \\
        i_* i^* N [-1] \ar[r] &
        j_! j^* N \ar[r] &
        N \ar[r] &
        i_* i^* N  
\\}}
\] 
According to Theorem~\ref{1C}~(a2), they provide weight
filtrations of $j_! j^* M$ and of $j_! j^* N$, respectively. By orthogonality
(condition~(3) in Definition~\ref{1A}), any morphism
from $i_* i^* M [-1]$ to $N$ is zero. Therefore, the above 
diagram can be completed to give a morphism of exact triangles. 
\end{Proofof}

\begin{Rem}
Following the lines of part (a) of the above proof, one can show that
there is in fact a canonical bijection
between the isomorphism classes of extensions of $M_U$ to $X$
as Chow motives on the one hand, and isomorphism classes of weight
filtrations of the restriction of a cone of $j_! M_U \to j_* M_U$
to the complement $X - U$ on the other hand. 
\end{Rem}

Let us note a consequence of Theorem~\ref{1F}, which we
think of as useful even though it will not be used in the rest
of this paper.

\begin{Cor} 
Let $j: U \into X$ be an open immersion of schemes.
Let $N_U^1, N_U^2 \in \CHUM$ and $M^1, M^2 \in \DBXM$. Then the image of the inverse image
\[
j^* : \Hom_X \bigl( M^1 \otimes_X N^1 , M^2 \otimes_X N^2 \bigr) \longto 
\Hom_U \bigl( j^* M^1 \otimes_U N_U^1 , j^* M^2 \otimes_U N_U^2 \bigr)
\]
is independent of the extensions of $N_U^n$ to Chow motives $N^n$ over $X$, $n = 1 , 2$.
\end{Cor}

\begin{Proof}
Let $N_r^n \in \CHXM$, $r = 1,2$ be two extensions of $N_U^n$, $n = 1 , 2$. By Theorem~\ref{1F}~(b),
there are morphisms $\beta_1^n: N_1^n \to N_2^n$ and $\beta_2^n: N_2^n \to N_1^n$ extending
$\id_{N_U^n}$. But then,
\[
j^* : \Hom_X \bigl( M^1 \otimes_X N_1^1 , M^2 \otimes_X N_1^2 \bigr) \longto 
\Hom_U \bigl( j^* M^1 \otimes_U N_U^1 , j^* M^2 \otimes_U N_U^2 \bigr)
\]
factors through $\Hom_X ( M^1 \otimes_X N_2^1 , M^2 \otimes_X N_2^2 )$, and
\[
j^* : \Hom_X \bigl( M^1 \otimes_X N_2^1 , M^2 \otimes_X N_2^2 \bigr) \longto 
\Hom_U \bigl( j^* M^1 \otimes_U N_U^1 , j^* M^2 \otimes_U N_U^2 \bigr)
\]
factors through $\Hom_X ( M^1 \otimes_X N_1^1 , M^2 \otimes_X N_1^2 )$.
\end{Proof}

\begin{Cor} \label{1Bei}
Let $j: U \into X$ be an open immersion of schemes.
Let $N_U \in \CHUM$ and $(p,q) \in \BZ^2$. Then the image of the inverse image
\[
j^* : \Hom_X \bigl( \one_X , N(p)[q] \bigr) \longto \Hom_U \bigl( \one_U , N_U(p)[q] \bigr)
\]
is independent of the extension of $N_U$ to a Chow motive $N$ over $X$.
\end{Cor}

\begin{Rem} \label{1Sch}
Corollary~\ref{1Bei} should be compared to Scholl's construction of ``integral'' sub-spaces of
motivic cohomology for Chow motives over local and global fields \cite[Sect.~1]{Sl}.
In fact, \emph{continuity} \cite[Thm.~14.2.5]{CD} implies that 
both statements of Theorem~\ref{1F} continue to hold when passing 
to the limit over all open sub-schemes of a given scheme $X$.
In particular, for any Dedekind domain $A$ with fraction field $K$, 
the restriction from $CHM(\Spec A)_\BQ$ to $CHM(\Spec K)_\BQ$
is essentially surjective and full.
This yields the cate\-gorial interpretation of \cite[Sect.~1]{Sl}.
It also
shows that Scholl's construction generalizes to the inclusion of a generic point
of \emph{any} scheme $X$ (always in the sense of our conventions), 
which may thus be chosen differently from the spectrum of a Dedekind domain.
\end{Rem}

We finish this section with a discussion of the behaviour of weights
under duality. Fix $X$, 
and suppose that the structure morphism of $X$
factors over a morphism $a : X \to S'$ with regular target. 
Fix an integer $d$, put
\[
R := a^! \one_{S'}(-d)[-2d] \in \DBcXM \; ,
\]
and form the local duality functor $\BD_X$ with respect to this choice of $R$
(see point~(D) above). Part~(a) of the following is contained in \cite[Cor.~3.9]{H};
statements~(a)--(c) are proved for regular $X$ in
\cite[Cor.~2.2.5]{H2}.

\begin{Thm} \label{1G}
Let $n$ be an integer, and consider the functor
\[
\BD_X : \DBcXM^{opp} \longto \DBcXM \; .
\]
(a)~$\BD_X$ maps $\DBcXM^{opp}_{w \le n}$ to $\DBcXM_{w \ge -n}$. \\[0.1cm]
(b)~$\BD_X$ maps $\DBcXM^{opp}_{w \ge n}$ to $\DBcXM_{w \le -n}$. \\[0.1cm]
(c)~$\BD_X$ maps $\CHXM^{opp}$ to $\CHXM$. 
\end{Thm}

Given that $\id_X = \BD_X^2$ on $\DBcXM$, we see that $\BD_X$ actually 
induces equivalences of categories
$\DBcXM^{opp}_{w \le n} \cong \DBcXM_{w \ge -n}$ \emph{etc.}

\medskip

\begin{Proofof}{Theorem~\ref{1G}}
The thick triangulated category $\DBcXM$ is generated by its heart $\CHXM$, 
and $\BD_X$ inverts the sign of the shifts. Therefore, it suffices 
to prove part~(c). By \cite[Thm.~3.3~(ii)]{H} (see Remark~\ref{1E}), it
is enough to prove that for any proper morphism 
$f: Y \to X$ with regular source $Y$, and any integer $p$, the 
constructible Beilinson motive
\[
\BD_X \bigl( f_! \one_Y (p)[2p] \bigl) 
\]
is actually a Chow motive. From the formulae recalled in 
point~(D) above, 
\[
\BD_X \bigl( f_! \one_Y (p)[2p] \bigl) \cong
f_! \BD_Y \bigl( \one_Y (p)[2p] \bigl) 
\]
(recall that $f$ is proper), 
provided $\BD_Y$ is formed with respect to $f^! R$.  
But $f^! R = (a\circ f)^! \one_{S'}(-d)[-2d]$, hence
\[
\BD_Y \bigl( \one_Y (p)[2p] \bigl) \cong
(a\circ f)^! \one_{S'}(-(d+p))[-2(d+p)] \; .
\]
$Y$ has a finite Zariski covering by connected
quasi-projective schemes $Y_i$ over $S'$.
Therefore (still thanks to point~(D) above),
the restriction to any $Y_i$ of
the motive $(a\circ f)^! \one_{S'}(-(d+p))[-2(d+p)]$ is
isomorphic to $\one_{Y_i}(m)[2m]$, for some integer $m$.
In particular, we see that $\BD_Y \bigl( \one_Y (p)[2p] \bigl)$
is Zariski-locally of weight zero. The two localization
triangles, together
with the $w$-exactness properties from Theorem~\ref{1C}~(a2)
then show that $\BD_Y \bigl( \one_Y (p)[2p] \bigl)$ is 
itself of weight zero.
Again by Theorem~\ref{1C}~(a2), the same is true for
its image under $f_!$.
\end{Proofof}

%%% Local Variables:
%%% mode: latex
%%% TeX-master: "head"
%%% End:

\bigskip
%\include{Sec2}
%%%%%%%%%%%%%%%%%%%%%%%%%%%%%%%%%%%%%%%%%%%%%%%%%%%%%%%%%%%%%%%%%%%%%%%
%
%  Section 2
%
%%%%%%%%%%%%%%%%%%%%%%%%%%%%%%%%%%%%%%%%%%%%%%%%%%%%%%%%%%%%%%%%%%%%%%%

\section{Definition of the motivic intersection complex}
\label{2}

%%%%%%%%%%%%%%%%%%%%%%%%%%%%%%%%

%%%%%%%%%%%%%%%%%%%%%%%%%%%%%%%%

Fix a scheme $X$. Since (by the conventions fixed in the beginning)
$X$ is excellent, there is an open immersion $j: U \into X$
whose image $U$ is dense in $X$, and regular. 
Recall that by Theorem~\ref{1C}~(a1),
the Beilinson motive $\one_U$ belongs to $\CHUM$, and that by
Theorem~\ref{1F}~(a), it can be extended to $\CHXM$.

\begin{Def} \label{2A}
A pair $(j_{!*} \one_U , \alpha)$ is called \emph{motivic intersection
complex} on $X$ if the following conditions are satisfied.
\begin{enumerate}
\item[(1)] The object $j_{!*} \one_U$ belongs to $\CHXM$, and 
\[
\alpha : j^* j_{!*} \one_U \isoto \one_U
\]
is an isomorphism in $\CHUM$.
\item[(2)] The morphism induced by $\alpha$,
\[
j^* : \End_{\CHXM} \bigl( j_{!*} \one_U \bigr) \longto
\End_{\CHUM} \bigl( \one_U \bigr)
\]
is injective.
\end{enumerate}
\end{Def}

Given that $j^*$ is full (Theorem~\ref{1F}), axiom~(2) is equivalent
to requiring the restriction from $\End_{\CHXM} ( j_{!*} \one_U )$ to 
$\End_{\CHUM} ( \one_U )$ to be bijective. \\

Denote by $i$ the closed immersion of the complement $Z$ (with the
reduced structure, say) into $X$.

\begin{Rem} \label{2rem}
When $S = \Spec k$ for a finite field $k$ of characteristic $p$,
let us consider the formalism of weights on \emph{perverse $\ell$-adic sheaves},
for $\ell \ne p$ \cite[Sect.~5]{BBD}. \\[0.1cm]
(a)~One of the main results
from \loccit \ states that $j_{!*}$ is a functor which
transforms perverse sheaves which are pure of a given weight
into perverse sheaves which are pure of the same weight
\cite[Cor.~5.4.3]{BBD}. In particular, the intersection complex 
$j_{!*} \BQ_\ell$ is indeed pure of weight $0$. \\[0.1cm]
(b)~Localization implies that the kernel of  
\[
j^* : \End_X \bigl( j_{!*} \BQ_\ell \bigr) \longto
\End_U \bigl( \BQ_\ell \bigr)
\]
is a quotient of the group
\[
\Hom_Z \bigl( i^* j_{!*} \BQ_\ell , i^! j_{!*} \BQ_\ell \bigr) \; .
\]
But this group is zero since with respect to the 
\emph{perverse $t$-structure},
the object $i^* j_{!*} \BQ_\ell$ in concentrated in degrees $\le -1$,
while $i^! j_{!*} \BQ_\ell$ is in degrees $\ge 1$ \cite[Cor.~1.4.25]{BBD}.
\end{Rem}

Since the (perverse) $t$-structure for Beilinson motives is not known
to exist in general, the na\"{\i}ve generalization of the definition of the
intersection complex is not possible
(but see \cite[Sect.~3]{S} for the case of Artin--Tate
motives over a number ring). Definition~\ref{2A} circumvents
this problem by replacing the use of a $t$-structure by the use
of the motivic weight structure! 

\begin{Ex} \label{2ex}
If $X$ is regular, then $(\one_X , \id)$ is a motivic
intersection complex, as follows from the relation to $K$-theory
(see Section~\ref{1}, point~(A)), and from the invariance under
passage from $X$ to its reduced structure
$X_{red}$ (see Section~\ref{1}, point~(E)). 
Indeed, the restriction
\[
j^* : \End_{\CHXM} \bigl( \one_X \bigr) \longto
\End_{\CHUM} \bigl( \one_U \bigr)
\]
then corresponds to 
\[
j^* : \Gr_\gamma^0 K_0(X_{red})_\BQ \longto \Gr_\gamma^0 K_0(U_{red})_\BQ \; .
\]
The latter is an isomorphism since both sides
are canonically isomorphic (\emph{via} the rank)
to $r$ copies of $\BQ$, where $r$ is the number of
connected components of $X$, which coincides with the 
number of connected components of $U$ (recall that $U$ is dense in $X$). 
\end{Ex}

The same argument shows the following.

\begin{Prop} \label{2B}
The motivic intersection complex does not depend on the choice
of dense open regular sub-scheme of $X$. More precisely,
if $V$ is a dense open regular
sub-scheme of $X$ contained in $U$, 
and if $(j_{!*} \one_U , \alpha)$
is a motivic intersection complex with respect to $U$,
then 
\[
\bigl( j_{!*} \one_U , \alpha_{\tei V} \bigr)
\]
is a motivic intersection complex with respect to $V$. 
\end{Prop}

The proof of the following requires more efforts.

\begin{Prop} \label{2C}
The motivic intersection complex is compatible with 
restriction to open sub-schemes $W$ of $X$. More precisely,
if $(j_{!*} \one_U , \alpha)$
is a motivic intersection complex on $X$, then 
\[
\bigl( (j_{!*} \one_U)_{\tei W} , \alpha_{\tei W \cap U} \bigr)
\]
is a motivic intersection complex on $W$. 
\end{Prop} 

\begin{Proof}
Assume first that $W$ is dense in $X$. 
By Proposition~\ref{2B}, we may assume $U$ to be contained in $W$.
Let 
\[
\beta_W : (j_{!*} \one_U)_{\tei W} \longto (j_{!*} \one_U)_{\tei W}
\]
be an endomorphism restricting trivially to $U$. 
The inverse image from $X$ to $W$ is full (Theorem~\ref{1F}~(b)),
therefore $\beta_W$ is the restriction to $W$ of an endomorphism
$\beta$ of $j_{!*} \one_U$.
By assumption, we have $j^* \beta = 0$.
Condition~(2) of Definition~\ref{2A} implies that $\beta = 0$.
Hence $\beta_W = 0$.

In the general case, we follow an argument due to G.~Ancona \cite{Anc}. First,
choose an open sub-scheme $W'$ of $X$, contained in the
complement of $W$, and such that $W \coprod W'$ is dense in $X$.
By the above, the restriction to $W \coprod W'$ of $j_{!*} \one_U$
is a motivic intersection complex. We are thus reduced to the case
where $X = W \coprod W'$. We leave it to the reader to show that
the restrictions to $W$ and $W'$ of $j_{!*} \one_U$ are then
motivic intersection complexes on $W$ and $W'$, respectively.
\end{Proof}

%%% Local Variables:
%%% mode: latex
%%% TeX-master: "head"
%%% End:

\bigskip
%\include{Sec3}
%%%%%%%%%%%%%%%%%%%%%%%%%%%%%%%%%%%%%%%%%%%%%%%%%%%%%%%%%%%%%%%%%%%%%%%
%
%  Section 3
%
%%%%%%%%%%%%%%%%%%%%%%%%%%%%%%%%%%%%%%%%%%%%%%%%%%%%%%%%%%%%%%%%%%%%%%%

\section{Basic properties}
\label{3}

%%%%%%%%%%%%%%%%%%%%%%%%%%%%%%%%

%%%%%%%%%%%%%%%%%%%%%%%%%%%%%%%%

We keep the previous setting. Thus, $X$ is fixed scheme, and
$j: U \into X$ the immersion of a dense open regular sub-scheme.
The complementary immersion is denoted by $i: Z \into X$.

\begin{Thm} \label{3A}
(a)~The motivic intersection complex is unique up to uni\-que isomorphism. 
\\[0.1cm]
(b)~If the motivic intersection complex $(j_{!*} \one_U , \alpha)$
exists, and if 
\[
\beta : j^* M \isoto \one_U
\]
is an isomorphism in $\CHUM$, with $M \in \CHXM$,
then $j_{!*} \one_U$ is (in general non-canonically) 
a direct factor of $M$. More precisely, there
is an isomorphism
\[
M \isoto j_{!*} \one_U \oplus i_* L_Z
\]
restricting to $\alpha^{-1} \circ \beta$ on $U$,
with $L_Z \in \CHZM$.
\end{Thm}

\begin{Proof}
Recall that the inverse image $j^*$ is full on $CHM(\argdot)$
(Theorem~\ref{1F}~(b)). Therefore, there exist
morphisms of Chow motives
\[
\varphi: j_{!*} \one_U \longto M \quad \text{and} \quad
\psi: M \longto j_{!*} \one_U
\] 
extending $\beta^{-1} \circ \alpha$ and $\alpha^{-1} \circ \beta$,
respectively. Observe that the composition $\psi \circ \varphi$
restricts to the identity on $\one_U$. Injectivity of 
\[
j^* : \End_{\CHXM} \bigl( j_{!*} \one_U \bigr) \longto
\End_{\CHUM} \bigl( \one_U \bigr)
\]
therefore implies that $\psi \circ \varphi = \id_{j_{!*} \one_U}$.

Similarly, $\varphi \circ \psi = \id_M$ if $(M , \beta)$
is another choice of motivic intersection complex; note
that in this case, the relations $\psi \circ \varphi = \id_{j_{!*} \one_U}$
and $\varphi \circ \psi = \id_M$
hold for \emph{any} choices of $\varphi$, $\psi$, 
meaning that they are actually unique.

In the general case, $\varphi \circ \psi$ is an idempotent endomorphism of $M$.
Since its restriction to $U$ is the identity, 
localization (see Section~\ref{1}, point~(E)) shows that its kernel is
necessarily a Chow motive of the form $i_* L_Z$.
The Beilinson motive $L_Z \in \DBcM(Z)$ equals
both $i^*i_* L_Z$ and $i^!i_* L_Z$. By Theorem~\ref{1C}~(a2),
it is of weight zero, hence a Chow motive over $Z$.
\end{Proof}

\begin{Rem}
In the context of perverse $\ell$-adic sheaves over schemes of finite
type over a finite field, the analogue of Theorem~\ref{3A}~(b)
(concerning pure complexes $M$ of $\ell$-adic sheaves on $X$) 
is a consequence of the Decomposition Theorem \cite[Thm.~5.4.5]{BBD}.
As illustrated by our proof,
the formalism of weight structures yields a structural reason
for the non-canonicity of the isomorphism of \loccit .
\end{Rem}

According to Theorem~\ref{3A}~(b),
the motivic intersection complex (provided it exists) is indeed
minimal among all possible extensions of $\one_U$
to a Chow motive over $X$. Furthermore, our result suggests
a possible stra\-te\-gy for its construction: first, 
use Theorem~\ref{1F}~(a) to \emph{choose}
any extension $M \in \CHXM$ of $\one_U$;
then, \emph{choose} idempotent endomorphisms 
of $M$ to split off direct factors of the shape $i_* L_Z$,
until no such factor is left. 
Note that it is not clear that the result is independent
of the choices (of $M$ and of the splittings) made in this process. 
Nor is it clear that the result actually satisfies axiom~(2)
of Definition~\ref{2A}.
We plan to elaborate on this elsewhere. 

\begin{Rem} \label{3B}
Let $X = X_1 \cup X_2$ be a covering by two dense open sub-schemes,
and assume that
the motivic intersection complexes on $X_1$ and on $X_2$ exist.
Using Proposition~\ref{2C}, Theorem~\ref{3A} and Theorem~\ref{1F}, one can
show that they can be glued along $X_1 \cap X_2$ 
to give $(M , \alpha)$, with $M \in \CHXM$, and 
\[
\alpha : j^* M \isoto \one_U \; .
\]
\end{Rem}

\begin{Exo} \label{3pro}
In the situation of Remark~\ref{3B},
show that $(M , \alpha)$ satisfies axiom~(2) of Definition~\ref{2A}. 
\end{Exo}

There is one specific case where we know the solution
to Problem~\ref{3pro}. It is worthwhile to spell it out. 

\begin{Cor} \label{3loc}
Assume that $X$ is regular, and that $M \in \DBXM$ is Nisnevich-locally
isomorphic to $\one_X$, \emph{i.e.}, there is a finite Nisnevich covering
of $X$ by schemes $U_n$ such that $M_{\tei U_n} \cong \one_{U_n}$
for all $n$. Then $X \cong \one_X$. 
\end{Cor}

\begin{Proof}
The separation property of $\DBM(\argdot)$ \cite[Def.~2.1.11]{CD}
and the $w$-exactness properties from Theorem~\ref{1C}~(a2)
allow to control the weights of $M$ locally for the smooth topology;
in particular, our assumptions imply that $M \in \CHXM$.
The given covering of $X$ can be refined to construct a dense open sub-scheme
$j : U \into X$ and an isomorphism
\[
\beta : j^* M \isoto \one_U \; .
\]
By Example~\ref{2ex} and Theorem~\ref{3A}~(b),
\[
M \cong \one_X \oplus i_* L_Z
\]
for some $L_Z \in \CHZM$. On each $U_n$, the restriction
\[
\one_{U_n} \oplus (i_* L_Z)_{\tei U_n}
\]
has the same endomorphisms as $\one_{U_n}$. Therefore,
$(i_* L_Z)_{\tei U_n} = 0$. Separation implies that $i_* L_Z = 0$.
\end{Proof}

Note that if $M$ is Zariski-locally isomorphic to $\one_X$, then
separation can be replaced by an application of the two localization triangles.

\begin{Rem}
Corollary~\ref{3loc} and \cite[Prop.~14.3.29~(i)]{CD} can be employed
to show that on a regular scheme $X$, two $\otimes$-invertible
objects of $\DBcXM$ are isomorphic as soon as they are Nisnevich-locally isomorphic. 
\end{Rem}

\begin{Cor}[Absolute purity] \label{3abs}
Let $a : X \to S'$ be a morphism
of pure relative dimension 
$d$: for any irreducible component of $S'$ of dimension $n$, 
its pre-image under $a$ is of pure dimension $n+d$.
Assume that both $X$ and $S'$ are regular.
Then there is an isomorphism
\[
a^! \one_{S'}(-d)[-2d] \cong \one_X \; . 
\]  
\end{Cor}

\begin{Proof}
Cover $X$ by open sub-schemes which are quasi-projective over $S'$.
The discussion from point~(D) of Section~\ref{1}
then shows that the assumption of Corollary~\ref{3loc}
is satisfied (even Zariski-locally) for $M = a^! \one_{S'}(-d)[-2d]$.
\end{Proof}

Let us come back to the general situaton, \emph{i.e.},
drop the regularity assumption on $X$.
We aim at a motivic analogue of
\cite[Prop.~2.1.17]{BBD} which states that $j_{!*} \CF$
is auto-dual on $X$ provided that $\CF$ is auto-dual on $U$.
In order to have the motivic analogue of that assumption satisfied
for $\one_U$, we suppose that the structure morphism of $X$
factors over a morphism $a : X \to S'$ with regular target. We also suppose
that $a$ is of pure relative dimension $d$. Put
\[
R := a^! \one_{S'}(-d)[-2d] \in \DBcXM \; ,
\]
and form the local duality functor $\BD_X$ with respect to this choice of $R$.
By absolute purity (Corollary~\ref{3abs}),
there is an isomorphism
\[
\gamma : \one_U \isoto j^* R \; .
\]  
We thus have
\[
\BD_X = \uHom_X (\argdot , R) \; ,
\]
and composition with $\gamma$ is an isomorphism of functors 
\[
\uHom_U (\argdot , \one_U) 
\isoto \uHom_U (\argdot , j^* R) = \BD_U \; .
\]
When evaluated on $\one_U$, this gives
\[
\gamma_*: \one_U = \uHom_U (\one_U , \one_U) \isoto \BD_U(\one_U) \; ;
\]
it is in this precise sense that $\one_U$ is auto-dual. 
Theorem~\ref{3A} has the following formal consequence.

\begin{Cor}[Auto-duality] \label{3C}
If the motivic intersection complex 
exists, then it is auto-dual. More precisely, there is a unique
isomorphism
\[
j_{!*} \one_U \isoto \BD_X(j_{!*} \one_U)  
\]
compatible with $\alpha$ and $\gamma_*$ in the sense that
its restriction to $U$ equals the composition
\[
\BD(\alpha) \circ \gamma_* \circ \alpha : 
j^* j_{!*} \one_U \isoto \BD_U(j^* j_{!*} \one_U) \; .
\]
\end{Cor}

\begin{Proof}
By Theorem~\ref{1G}~(c), $\BD_X(j_{!*} \one_U)$ is a Chow motive over $X$.
Define 
\[
\beta := \gamma_*^{-1} \circ \BD(\alpha^{-1}) : 
j^* \BD_X(j_{!*} \one_U) = \BD_U(j^* j_{!*} \one_U) \isoto \one_U \; .
\]
With this choice, axiom~(2) of Definition~\ref{2A} is satisfied; 
indeed, by adjunction,
\[
\End_{\CHXM} \bigl( \BD_X(j_{!*} \one_U) \bigr) =
\Hom_{\CHXM} \bigl( j_{!*} \one_U , \BD_X^2(j_{!*} \one_U) \bigr) \; ,
\]
and $\id_X = \BD_X^2$ on $\CHXM$. 
Our claim follows from Theorem~\ref{3A}~(a).
\end{Proof}

\begin{Def} 
Assume that the motivic intersection complex $(j_{!*} \one_U , \alpha)$
exists. The pairing
\[
j_{!*} \one_U \otimes_X j_{!*} \one_U \longto a^! \one_{S'}(-d)[-2d]
\]
obtained by adjunction from the auto-duality isomorphism
is called the \emph{motivic intersection pairing}.
\end{Def}

By definition, the motivic intersection pairing is non-degenerate
in the sense that its adjoint is an isomorphism. Applying $a_!$
to the first component of its source, and $a_*$ to the second,
we get
\[
a_! j_{!*} \one_U \otimes_{S'} a_*j_{!*} \one_U \; ,
\]
which maps (isomorphically, by the projection formula 
\cite[Thm.~2.4.21~(v)]{CD}) to 
\[
a_! \bigl( j_{!*} \one_U \otimes_X a^* a_*j_{!*} \one_U \bigr) \; ,
\]
and finally, \emph{via} the adjunction $(a^*,a_*)$, to
\[
a_! \bigl( j_{!*} \one_U \otimes_X j_{!*} \one_U \bigr) \; .
\] 
Composition with $a_!$ of the intersection pairing,
and application of the adjunction $(a_!,a^!)$ yields the pairing
\[
a_! j_{!*} \one_U \otimes_{S'} a_*j_{!*} \one_U \longto \one_{S'}(-d)[-2d] \; .
\]
It is non-degenerate since by construction, its adjoint is the isomorphism
\[
a_! j_{!*} \one_U \isoto \BD_{S'} a_*j_{!*} \one_U
\]
obtained from $a_!$ of auto-duality and the formula $\BD_{S'} a_* = a_! \BD_X$
(see Section~\ref{1}, point~(D)). In particular, we get the motivic analogue
of Poincar\'e duality for intersection cohomology.

\begin{Cor} \label{3D}
Assume that the motivic intersection complex $(j_{!*} \one_U , \alpha)$
exists, and that the morphism $a: X \to S'$ is proper. 
Then $a_! j_{!*} \one_U$ is auto-dual.
\end{Cor}

Note that under the assumptions of Corollary~\ref{3D}, 
the object $a_! j_{!*} \one_U$
is a Chow motive over $S'$ (Theorem~\ref{1C}~(a2)). \\

Here are the few cases where we actually 
know the hypotheses of Theorem~\ref{3A}~(b) and Corollaries~\ref{3C}
and \ref{3D} to be satisfied.

\begin{Thm} \label{3Z}
The motivic intersection complex exists in the following cases. \\[0.1cm]
(a)~The normalization $\nX$ of the reduced scheme underlying $X$ is regular. \\[0.1cm]
(b)~$X$ is of dimension at most two, and the residue fields of the singular points
of $\nX$ are perfect.  
\end{Thm} 

The proof of Theorem~\ref{3Z} will be given in the next section.

\begin{Rem} \label{3E}
Let us discuss the case $S = \Spec \BC \,$. We consider the Hodge theoretic realization
\[
R : \DBcM(\Spec \BC) \longto D^b \bigl( \MHS_\BQ \bigr)
\]
(\cite[Sect.~2.3 and Corrigendum]{Hu}; see \cite[Sect.~1.5]{DG}
for a simplification of this approach). 
Here, $D^b ( \MHS_\BQ )$ is the bounded derived category 
of the Abelian category $\MHS_\BQ$ of mixed graded-polarizable $\BQ$-Hodge structures.
It is reasonable to expect the Hodge realization to extend
to the relative setting, yielding exact, monoidal functors
\[
R : \DBcM(X) \longto D^b \bigl( \MHM_\BQ X \bigr)
\]
for all schemes $X$ over $\BC$. Here, $D^b ( \MHM_\BQ X )$ is the bounded derived category 
of algebraic mixed $\BQ$-Hodge modules on $X$ \cite{Sa}. Let us \emph{assume} such an 
extension $R$ to exist, and to be compatible
with the six operations from \cite{CD} and from \cite{Sa}. 
According to \cite[Prop.~2.7~I]{Bo3}, the category $D^b ( \MHM_\BQ X )$
carries a weight structure, with $D^b ( \MHM_\BQ X )_{w \le 0}$ and $D^b ( \MHM_\BQ X )_{w \ge 0}$
equal to the sub-categories of complexes of
Hodge modules of weights $\le 0$ and $\ge 0$
(in the sense of \cite[Def.~4.5]{Sa}), respectively. The Hodge theoretic realization
is then necessarily $w$-exact: indeed, since $\DBcM(X)$
is generated by its heart (Remark~\ref{1E}), it suffices to show that Chow motives over $X$
are mapped to Hodge modules which are pure of weight zero. This in turn
follows from the explicit 
description of $\CHXM$, and from the $w$-exactness properties
of the six operations on algebraic Hodge modules: any Chow motive over $X$ 
is a direct factor of 
\[
f_! \one_Y (p)[2p] \; ,
\]
for a proper morphism $f: Y \to X$ with regular source $Y$, and an integer $p$
(Remark~\ref{1E}). Its image under $R$ is therefore a direct factor of 
\[
f_! \BQ_Y^H (p)[2p] \; .
\]
The $\BC$-scheme $Y$ is
regular, hence smooth over $\BC$. By \cite[Thm.~3.27]{Sa}, the variation
of Hodge structure $\BQ_Y^H(p)$ on $Y$ is a complex of algebraic Hodge modules;
as such, it is pure of weight $-2p$. Therefore, its shift $\BQ_Y^H(p)[2p]$ is pure of weight zero.
But by \cite[Sect.~(4.5.2)]{Sa}, $f_! = f_*$ is $w$-exact. 

The canonical $t$-structure on $D^b ( \MHM_\BQ X )$ allows to define the
Hodge theoretic intersection complex $IC_X \BQ^H$ on $X$ \cite[Sect.~4.5]{Sa}.
Due to the normalization we chose for the motivic intersection complex, we define 
\[
j_{!*} \BQ_U^H := IC_X \BQ^H [-d]
\]
if $X$ is of pure dimension $d$. (Thus, 
$j_{!*} \BQ_U^H$ is a complex of Hodge modules concentrated in degree $d$.) 
According to \cite[Sect.~4.5]{Sa}, $j_{!*} \BQ_U^H$ is pure of weight zero,
and extends $\BQ_U^H$. 
It satisfies the
Hodge theoretic analogue of axiom~(2) of Definition~\ref{2A}. From the 
Hodge theoretic analogue of Theorem~\ref{3A}~(b), we conclude that
the realization $R (j_{!*} \one_U)$ of the motivic intersection complex
contains $j_{!*} \BQ_U^H$ as a direct factor.

When $S = \Spec k$ is the spectrum of a finite field $k$, similar remarks
apply to perverse $\ell$-adic sheaves over $k$-schemes.
\end{Rem}

\begin{Exo}
In the situation of Remark~\ref{3E}, show the equality 
\[
R (j_{!*} \one_U) =j_{!*} \BQ_U^H \; .
\]
Note that it implies that the intersection motive of $X$ realizes to give
(the complex computing) intersection cohomology of $X$.
\end{Exo}

%%% Local Variables:
%%% mode: latex
%%% TeX-master: "head"
%%% End:

\bigskip
%\include{Sec4}
%%%%%%%%%%%%%%%%%%%%%%%%%%%%%%%%%%%%%%%%%%%%%%%%%%%%%%%%%%%%%%%%%%%%%%%
%
%  Section 4
%
%%%%%%%%%%%%%%%%%%%%%%%%%%%%%%%%%%%%%%%%%%%%%%%%%%%%%%%%%%%%%%%%%%%%%%%

\section{On the problem of existence}
\label{4}

%%%%%%%%%%%%%%%%%%%%%%%%%%%%%%%%

%%%%%%%%%%%%%%%%%%%%%%%%%%%%%%%%

We keep the situation considered before: $X$ is a scheme, and
$j: U \into X$ is the immersion of a dense open regular sub-scheme.
The complementary immersion is denoted by $i: Z \into X$.

\begin{Prop} \label{4B}
Let $Y$ be a scheme, $h: V \into Y$ an open immersion,
and $k: T \into Y$ the complement. Let $L$ and $N$
be objects of $\DBYM$. Assume that
$\Hom_T ( k^* L , k^! N ) = 0$. Then the restriction
\[
h^* : \Hom_Y \bigl( L , N \bigr) \longto
\Hom_V \bigl( h^* L , h^* N \bigr) 
\]
is injective. 
\end{Prop}

\begin{Proof}
Either one of the
localization triangles implies that the kernel of  
\[
h^* : \Hom_Y \bigl( L , N \bigr) \longto
\Hom_V \bigl( h^* L , h^* N \bigr) 
\]
is a quotient of $\Hom_T ( k^* L , k^! N )$.
\end{Proof}

In the setting of interest for us, Proposition~\ref{4B}
implies the following. 

\begin{Cor} \label{4A}
Assume that $M \in \CHXM$ is given, together with an isomorphism
\[
\alpha : j^* M \isoto \one_U \; .
\]
If 
\[
\Hom_Z \bigl( i^* M , i^! M \bigr) = 0 \; ,
\]
then $(M,\alpha)$ equals the motivic intersection complex on $X$.
\end{Cor}

Replacing axiom~(2) of Definition~\ref{2A} by the vanishing of  
\[
\Hom_Z \bigl( i^* j_{!*} \one_U , i^! j_{!*} \one_U \bigr)
\]
might possibly provide a ``better'' definition of the motivic intersection
complex. At least, the proof of Theorem~\ref{3Z} will consist in
showing this vanishing.
In order to do so, 
the following principle will be frequently used.

\begin{Cor} \label{4C}
Let $Y$ be a scheme, $h: V \into Y$ an open immersion,
and $k: T \into Y$ the complement. Let $L$ and $N$
be objects of $\DBYM$. 
If 
\[
\Hom_V \bigl( h^* L , h^* N \bigr) \quad \text{and} \quad
\Hom_T \bigl( k^* L , k^! N \bigr) = 0 \; ,
\]
then
\[
\Hom_Y \bigl( L , N \bigr) = 0 \; .
\]
\end{Cor}

Successive applications of this principle show that the vanishing
assumption of Corollary~\ref{4A} can be verified on a finite stratification.

\begin{Ex}
We get another proof of the equality 
``$j_{!*} \one_U = \one_X$'' for regular $X$ (Example~\ref{2ex}):
choose a stratification of $Z$ by regular sub-schemes $T$.
Then apply absolute purity and the relation to $K$-theory to see that 
\[
\Hom_T \bigl( i^* \one_X , i^! \one_X \bigr) = 0 
\]
for each $T$. 
\end{Ex}

Let us turn to the proof of Theorem~\ref{3Z}. We may assume
that $X$ is reduced. Denote by 
\[
p : \nX \longto X
\]
the normalization of $X$; note that $p$ is finite since $X$
is excellent. Note also that $j$ factors uniquely 
through an open immersion $\nj$ into $\nX$,
identifying $U$ with its pre-image under $p$.
Part~(a) of Theorem~\ref{3Z} is contained
in the following.

\begin{Prop} \label{4D}
Assume that $\nX$ is regular. Then $(p_! \one_{\nX},\id)$
equals the motivic intersection complex on $X$.
\end{Prop} 

\begin{Proof}
First, note that $\nX$ being supposed regular,
the Beilinson motive $\one_{\nX}$ is indeed a Chow motive
(Theorem~\ref{1C}~(a1)). Since $p_! = p_*$, 
the same is true for $p_! \one_{\nX}$
(Theorem~\ref{1C}~(a2)).
Define $\nnX:= \nX \times_X \nX$, denote by $p_1$ and $p_2$
the projections of $\nnX$ to the two factors $\nX$, and by
$P$ the projection to $X$. Base change (see Section~\ref{1}, point~(B))
and adjunction,
first from $Z$ to $p^{-1}(Z)$ and then to $P^{-1}(Z)$, show that
\[
\Hom_Z \bigl( i^* p_! \one_{\nX} , i^! p_! \one_{\nX} \bigr) = 
\Hom_{P^{-1}(Z)} 
\bigl( i^{nn,*} p_2^* \one_{\nX} , i^{nn,!} p_1^! \one_{\nX} \bigr) \; ,
\]
where we let $\nni$ denote the immersion of $P^{-1}(Z)$ into $\nnX$.
Let $k : T \into P^{-1}(Z)$ be a regular connected locally 
closed sub-scheme.
It is necessarily quasi-finite over $Z$.
In particular, its relative dimension $e$ over $\nX$ via $p_1$
is strictly negative. As recalled in Section~\ref{1}, point~(D),
\[
k^! i^{nn,!} p_1^! \one_{\nX} = \one_T(e)[2e] \; .
\]
Of course,
\[
k^* i^{nn,*} p_2^* \one_{\nX} = \one_T \; .
\]
The relation to $K$-theory shows that
\[
\Hom_T 
\bigl( \one_T , \one_T(e)[2e] \bigr) = 0.
\]
Indeed, the graded object $\Gr_\gamma^e K_0(T)$ is zero since the
gamma filtration is concentrated in non-negative degrees.
Now apply Corollaries~\ref{4C} and \ref{4A}.
\end{Proof}

\begin{Exo} 
Without the regularity assumption on $\nX$,
show that  
\[
j_{!*} \one_U = p_! j_{!*}^{norm} \one_U \; ,
\]
whenever the motivic intersection complex
$j_{!*}^{norm} \one_U$ on $\nX$ exists. 
\end{Exo}

In order to prove part~(b) of Theorem~\ref{3Z}, note first that for reduced schemes
$X$ of dimension at most one, the normalization $\nX$ is regular.
For the rest of this section, let us therefore assume that $X$ is a reduced surface (\emph{i.e.},
all irreductible components of $X$ are integral and of dimension two),
and that the residue fields of the singular points of $\nX$ are perfect. \\

Let us start by the construction of $j_{!*}^{norm} \one_U$ on $\nX$. 
It is a variant of the construction from \cite{CM} for surfaces defined over a field.
By Proposition~\ref{2B}, we may perform the computation after replacing 
$U$ by the regular locus $V$ of $\nX$. Since $\nX$ is regular in codimension one,
the complement $Z'$ of $V$ (with the reduced structure) is finite; in fact, by our assumption, 
$Z'$ is the spectrum of a finite product of perfect fields. By
Abhyankar's result on resolution of singularities in dimension two 
\cite[Theorem]{L2}, $X$ can
be desingularized. In addition 
(see the discussion in \cite[pp.~191--194]{L1}), by further blowing
up possible singularities of (the components of) the pre-image $D$ of $Z'$, 
it can be assumed to be a divisor with
normal crossings, whose irreducible components are regular.
Fix such a resolution, that is, fix the following diagram, assumed
to be cartesian:
\[
\vcenter{\xymatrix@R-10pt{
        V \ar@{^{ (}->}[r] \ar@{=}[d] &
        \Xp \ar@{<-^{ )}}[r]^-{\ip} \ar[d]_-\pi &
        D \ar[d]^\pi \\
        V \ar@{^{ (}->}[r] &
        \nX \ar@{<-^{ )}}[r]^-{i'} &
        Z'
\\}}
\]
where $\pi$ is proper (and birational), $\Xp$ is regular, and $D$ is a divisor with
normal crossings, whose irreducible components $D_m$ are regular. \\

Thus, the $D_m$ are regular curves over perfect fields (the points of $Z'$).
Therefore, they are smooth. In addition, they are proper. Denote by $\ip_m$
the closed immersion of $D_m$ into $\Xp$, and by $\pi_m$
the restriction of $\pi$ to $D_m$. The classical
theory of Chow motives yields
canonical (split) sub-objects $(\pi_{m,!} \one_{D_m})^0$ and (split)
quotients $(\pi_{m,!} \one_{D_m})^2$ of $\pi_{m,!} \one_{D_m}$. 
The adjunctions $\id_{\Xp} \to \ip_{m,*} \ip_m^*$ and $\ip_{m,*} \ip_m^! \to \id_{\Xp}$,
and absolute purity for $\ip_m$ yield canonical morphisms 
\[
\ip^* : \pi_! \one_{\Xp} \longto \bigoplus_m \pi_{m,!} \one_{D_m} 
\longonto \bigoplus_m \bigl( \pi_{m,!} \one_{D_m} \bigr)^2
\]
and
\[
\ip_* : \bigoplus_m \bigl( \pi_{m,!} \one_{D_m} \bigr)^0 (-1)[-2]
\longinto \bigoplus_m \pi_{m,!} \one_{D_m}(-1)[-2] \longto \pi_! \one_{\Xp} 
\]
of Chow motives over $\nX$.

\begin{Prop} \label{4E}
(i) The composition $\alpha := {\ip}^*{\ip}_*$ is an isomorphism. \\[0.1cm]
(ii) The composition $\varepsilon:= {\ip}_*\alpha^{-1}{\ip}^*$ is an idempotent on $\pi_! \one_{\Xp}$.
Hence so is the difference $\id_{\pi_! \one_{\Xp}}-\varepsilon$. \\[0.1cm]
(iii) The image $\imm \varepsilon$ is canonically isomorphic to $\oplus_m ( \pi_{m,!} \one_{D_m} )^2$.
\end{Prop}

\begin{Proof}
The proof is formally identical to the one of \cite[Thm.~2.2]{W2}.
Observe that the non-degeneracy of
the intersection pairing on the components of $D$ holds since the proof 
\cite[p.~6]{M} carries over to the general context
of normal surfaces.
\end{Proof}

Note that the $(\pi_{m,!} \one_{D_m})^2$ restrict trivially to $V$. Therefore, the
image $\imm (\id_{\pi_! \one_{\Xp}}-\varepsilon)$ restricts to give $\one_V$.
Part~(b) of Theorem~\ref{3Z} is contained in the following.

\begin{Prop} \label{4F}
Assume that $X$ is a surface, and that the residue fields of the singular points
of $\nX$ are perfect. \\[0.1cm]
(a)~The pair $(\imm (\id_{\pi_! \one_{\Xp}}-\varepsilon),\id)$ equals the motivic
intersection complex on $\nX$. \\[0.1cm]
(b)~The pair $(p_! \imm (\id_{\pi_! \one_{\Xp}}-\varepsilon),\id)$ equals the motivic
intersection complex on $X$.
\end{Prop}

\begin{Proof}
It suffices to prove part~(b). Write $M := \imm (\id_{\pi_! \one_{\Xp}}-\varepsilon)$.
As in the proof of Proposition~\ref{4D},
let $\nnX= \nX \times_X \nX$, denote by $p_1$ and $p_2$
the projections of $\nnX$ to the two factors $\nX$, and by
$P$ the projection to $X$. Base change and adjunction show that
\[
\Hom_Z \bigl( i^* p_! M , i^! p_! M \bigr) = 
\Hom_{P^{-1}(Z)} \bigl( i^{nn,*} p_2^* M , i^{nn,!} p_1^! M \bigr) \; ,
\]
where $\nni$ denotes the immersion of $P^{-1}(Z)$ into $\nnX$.
In order to apply Corol\-la\-ry~\ref{4A}, we need to show the vanishing of this group.
We shall repeatedly apply Corollary~\ref{4C}. In order to do so,
stratify $P^{-1}(Z)$ as follows: the open stratum (possibly empty) is the intersection
of the pre-images under $p_1$ and under $p_2$ of $V$ (which contains $U$), 
the closed stratum is the complement, which is a finite set of points.

If $k : T \into P^{-1}(Z)$ is a regular connected locally 
closed sub-scheme of the open stratum, then its relative dimension $e$ over $\nX$ via $p_1$
is strictly negative. Since $p_1(T)$ and $p_2(T)$ are contained in $V$, and
$\imm (\id_{\pi_! \one_{\Xp}}-p)$ restricts to $\one_V$, we argue as in the proof of
Proposition~\ref{4D} to see that   
\[
\Hom_T \bigl( k^*i^{nn,*} p_2^* M , k^!i^{nn,!} p_1^! M \bigr) = 0 \; .
\]
It remains to check the points $k : T \into P^{-1}(Z)$ of
the closed stratum. Depending on whether $p_2(T)$ is regular or not, we have
\[
k^*i^{nn,*} p_2^* M = \one_T 
\]
or
\[
k^*i^{nn,*} p_2^* M = \bigl( \pi_! \one_{D_{p_2(T)}} \bigr)^{\le 1} \; ,
\] 
where $D_{p_2(T)}$ is the base change of the exceptional divisor $D$ to $p_2(T)$, and
where the symbol $(\pi_! \one_{D_{p_2(T)}} \bigr)^{\le 1}$ denotes 
the kernel of the projection 
\[
\pi_! \one_{D_{p_2(T)}} \longonto 
(\pi_! \one_{D_{p_2(T)}} \bigr)^2 := \bigoplus_m (\pi_{m,!} \one_{D_{m,p_2(T)}})^2
\]
induced by $\varepsilon$.
Similarly (see Section~1, point~(D)), 
\[
k^!i^{nn,!} p_1^! M = \one_T(-2)[-4] 
\]
or
\[
k^!i^{nn,!} p_1^! M = \bigl( \pi_! (\ip_{\tei D_{p_1(T)}})^! \one_{\Xp} \bigr)^{\le 1} \; .
\] 
We thus need to distinguish four cases. As usual,
\[
\Hom_T \bigl( \one_T , \one_T(-2)[-4] \bigr) = 0 \; .
\]
In order to show that
\[
\Hom_T \bigl( \one_T , \bigl( \pi_! (\ip_{\tei D_{p_1(T)}})^! \one_{\Xp} \bigr)^{\le 1} \bigr) = 0 \; ,
\]
note that $\pi_! = \pi_*$, and that by adjunction,
\[
\Hom_T \bigl( \one_T , \pi_* (\ip_{\tei D_{p_1(T)}})^! \one_{\Xp} \bigr) 
= \Hom_{D_{p_1(T)}} \bigl( \one_{D_{p_1(T)}} , (\ip_{\tei D_{p_1(T)}})^! \one_{\Xp} \bigr) \; .
\]
Stratify $D_{p_1(T)}$ by regular sub-schemes, using that $\ip$ is of stricly positive
codimension, to see that the latter group is zero.

Similarly, 
\[
\Hom_T \bigl( \bigl( \pi_! \one_{D_{p_2(T)}} \bigr)^{\le 1} , \one_T(-2)[-4] \bigr) = 0 \; .
\]
It remains to consider
\[
\Hom_T \bigl( \bigl( \pi_! \one_{D_{p_2(T)}} \bigr)^{\le 1} , 
              \bigl( \pi_! (\ip_{\tei D_{p_1(T)}})^! \one_{\Xp} \bigr)^{\le 1} \bigr) \; .
\]
Adjunction yield an identification between 
\[
\Hom_T \bigl( \pi_! \one_{D_{p_2(T)}} , 
              \pi_! (\ip_{\tei D_{p_1(T)}})^! \one_{\Xp} \bigr) 
\]
and
\[
\Hom_{D_{p_2(T)} \times_T D_{p_1(T)}} \bigl( \one_{D_{p_2(T)} \times_T D_{p_1(T)}} , 
              pr_{D_{p_1(T)}}^! (\ip_{\tei D_{p_1(T)}})^! \one_{\Xp} \bigr) \; ,
\]
where $D_{p_2(T)} \times_T D_{p_1(T)}$ is the (singular) surface obtained by base
change over $T$ of the curves $D_{p_2(T)}$ and $D_{p_1(T)}$, and $pr_{D_{p_1(T)}}$
denotes the projection to $D_{p_1(T)}$. On regular sub-schemes contained
in the singular locus of $D_{p_2(T)} \times_T D_{p_1(T)}$, the same considerations
as before show that there are no non-zero morphisms of the required type.
Hence the assumption of Proposition~\ref{4B} is satisfied,
and 
\[
\Hom_{D_{p_2(T)} \times_T D_{p_1(T)}} \bigl( \one_{D_{p_2(T)} \times_T D_{p_1(T)}} , 
              pr_{D_{p_1(T)}}^! (\ip_{\tei D_{p_1(T)}})^! \one_{\Xp} \bigr) 
\]
injects into 
\[
\Hom_{{D_{p_2(T)}}_{reg} \times_T {D_{p_1(T)}}_{reg}} 
\bigl( \one_{{D_{p_2(T)}}_{reg} \times_T {D_{p_1(T)}}_{reg}} , 
pr_{{D_{p_1(T)}}_{reg}}^! (\ip_{\tei {D_{p_1(T)}}_{reg}})^! \one_{\Xp} \bigr) \; ,
\]
where the subscripts ${}_{reg}$ denote the regular loci. Both 
${D_{p_2(T)}}_{reg} \times_T {D_{p_1(T)}}_{reg}$ and $\Xp$ are regular surfaces, hence
\[
pr_{{D_{p_1(T)}}_{reg}}^! (\ip_{\tei {D_{p_1(T)}}_{reg}})^! \one_{\Xp} 
= \one_{{D_{p_2(T)}}_{reg} \times_T {D_{p_1(T)}}_{reg}} \; .
\]
Therefore,
\[
\Hom_{{D_{p_2(T)}}_{reg} \times_T {D_{p_1(T)}}_{reg}} 
\bigl( \one_{{D_{p_2(T)}}_{reg} \times_T {D_{p_1(T)}}_{reg}} , 
pr_{{D_{p_1(T)}}_{reg}}^! (\ip_{\tei {D_{p_1(T)}}_{reg}})^! \one_{\Xp} \bigr)
\]
equals $r$ copies of $\BQ$, where $r$ is the number of connected components
of ${D_{p_2(T)}}_{reg} \times_T {D_{p_1(T)}}_{reg}$. But the same result,
with compatible identifications is obtained
by computing 
\[
\Hom_T \bigl( \bigl( \pi_! \one_{D_{p_2(T)}} \bigr)^2 , 
              \bigl( \pi_! (\ip_{\tei D_{p_1(T)}})^! \one_{\Xp} \bigr)^2 \bigr) \; .
\]
\end{Proof}

%%% Local Variables:
%%% mode: latex
%%% TeX-master: "head"
%%% End:

\bigskip

%%%%%%%%%%%%%%%%%%%%%%%%%%%%%%%%%%%%%%%%%%%%%%%%%%%%%%%%%%%%%%%%%%%%%%%
%
%  Bibliography
%
%%%%%%%%%%%%%%%%%%%%%%%%%%%%%%%%%%%%%%%%%%%%%%%%%%%%%%%%%%%%%%%%%%%%%%%


\begin{thebibliography}{99}

\bibitem[An]{Anc}
G.~Ancona,
letter to the author dated Apr.~14, 2011.

\bibitem[Ay]{A}
J.~Ayoub,
{\it Les six op\'erations de Grothendieck et le formalisme des cycles
\'evanescents dans le monde motivique~(I--II)},
Ast\'erisque~{\bf 314--315},
Soc.\ Math.\ France (2007).

\bibitem[BBD]{BBD}
A.A.~Beilinson, J.~Bernstein, P.~Deligne,
{\it Faisceaux pervers},
in: B.~Teissier, J.L.~Verdier (eds.),
{\it Analyse et topologie sur les espaces singuliers (I)},
Ast\'erisque~{\bf 100},
Soc.\ Math.\ France (1982).

\bibitem[Bo1]{Bo1}
M.V.~Bondarko,
{\it Weight structures vs.\ $t$-structures; weight filtrations, 
spectral sequences, and complexes (for motives and in general)},
J.~$K$-Theory~{\bf 6} (2010), 387--504.

\bibitem[Bo2]{Bo2}
M.V.~Bondarko,
{\it Weight for relative motives; relation with mixed sheaves},
preprint, July~2010, version dated Sept.~21, 2010, 42~pages, 
available on arXiv.org under 
{\tt http://arxiv.org/abs/1007.4543} 

\bibitem[Bo3]{Bo3}
M.V.~Bondarko,
{\it Weights and $t$-structures: in general triangulated categories,
for $1$-motives, mixed motives, and for mixed Hodge complexes and modules},
preprint, November~2010, version dated Mar.~28, 2011, 27~pages, 
available on arXiv.org under 
{\tt http://arxiv.org/abs/1011.3507} 

\bibitem[CD]{CD}
D.-C.~Cisinski, F.~D\'eglise,
{\it Triangulated categories of mixed motives},
preprint, December~2009, 229~pages, 
available on arXiv.org under 
{\tt http://arxiv.org/abs/0912.2110} 

\bibitem[CM]{CM}
M.A.A.~de Cataldo, L.~Migliorini, 
{\it The Chow motive of semismall re\-solutions},
Math.\ Res.\ Lett.~{\bf 11} (2004), 151--170.

\bibitem[CH]{CH}
A.~Corti, M.~Hanamura,
{\it Motivic decomposition and intersection Chow groups. I},
Duke Math.\ J.~{\bf 103} (2000), 459--522.

\bibitem[DG]{DG}
P.~Deligne, A.B.~Goncharov,
{\it Groupes fondamentaux motiviques de Tate mixte},
Ann.\ Scient.\ ENS~{\bf 38} (2005), 1--56.

\bibitem[GHM]{GHM}
B.B.~Gordon, M.~Hanamura, J.P.~Murre, 
{\it Relative Chow--K\"unneth projectors for modular varieties},
J.\ reine angew.\ Math~{\bf 558} (2003), 1--14. 

\bibitem[H\'e1]{H}
D.~H\'ebert,
{\it Structures de poids \`a la Bondarko sur les motifs de Beilinson},
to appear in Compositio Math., 17~pages,
available on arXiv.org under 
{\tt http://arxiv.org/abs/1007.0219} 

\bibitem[H\'e2]{H2}
D.~H\'ebert,
{\it Complexe de Poids, Dualit\'e et Motifs de Beilinson},
preprint, October~2010, 17~pages, 
available on arXiv.org under 
{\tt http://arxiv.org/abs/1010.5469} 

\bibitem[Hu]{Hu}
A.~Huber,
{\it Realization of Voevodsky's motives},
J.\ of Alg.\ Geom.~{\bf 9} (2000), 755--799, 
{\it Corrigendum},
{\bf 13} (2004), 195--207.

\bibitem[L1]{L1}
J.~Lipman, {\it Introduction to resolution of singularities}, in R.~Hartshorne (ed.), {\it
Algebraic Geometry. Proceedings of the Symposium in Pure Mathematics of the AMS, held at Humboldt
State University, Arcata, California, July 29--August 16, 1974}, Proc.\ of Symp.\ in Pure
Math.~{\bf 29}, AMS (1975), 187--230.

\bibitem[L2]{L2}
J.~Lipman, {\it Desingularization of two-dimensional schemes}, Ann.\ of Math.~{\bf 107} (1978),
151--207.

\bibitem[M]{M}
D.~Mumford, {\it The topology of normal singularities of an algebraic surface and a criterion for
simplicity}, Publ.\ Math.\ IHES~{\bf 9} (1961), 5--22.

\bibitem[Sa]{Sa}
M.~Saito,
{\it Mixed Hodge Modules},
Publ.\ RIMS, Kyoto Univ.~{\bf 26} (1990), 221--333.

\bibitem[Sb]{S}
J.~Scholbach, 
{\it Mixed Artin--Tate motives over number rings},
preprint, March~2010, version dated Aug.~5, 2010, 28~pages, 
available on arXiv.org under 
{\tt http://arxiv.org/abs/1003.1267} 

\bibitem[S]{Sl}
A.J.~Scholl,
{\it Integral elements in $K$-theory and products of modular curves},
in B.B.~Gordon, J.D.~Lewis, S.~M\"uller-Stach, S.~Saito, N.~Yui (eds.),
{\it The Arithmetic and Geometry of Algebraic Cycles. 
Proceedings of the NATO Advanced Study Institute held as part of the 1998 CRM Summer School at Banff}, 
Kluwer Acad.\ Publ. (2000),  467--489.  

\bibitem[VSF]{VSF}
V.~Voe\-vodsky, A.~Suslin, E.M.~Friedlander,
{\it Cycles, Transfers, and Motivic Homology Theories},
Ann.\ of Math.\ Studies~{\bf 143}, Princeton Univ.\ Press 2000.

\bibitem[W1]{W1}
J.~Wildeshaus,
{\it On the interior motive of certain Shimura varieties: 
the case of Hilbert--Blumenthal varieties},
preprint, June~2009, version dated March~18, 2011, 31~pages, submitted, 
available on arXiv.org under 
{\tt http://arxiv.org/abs/0906.4239} 

\bibitem[W2]{W2}
J.~Wildeshaus,
{\it Pure motives, mixed motives and extensions of motives
associated to singular surfaces},
39 pages, to appear in M.~Kim, S.~Ramdorai, L.~Lafforgue, 
A.~Genestier, B.-C.~Ngo (eds.), 
{\it Autour des motifs. Actes de l'\'ecole d'\'et\'e franco-asiatique 
de g\'eom\'etrie alg\'ebrique et th\'eorie des nombres, IHES, Juillet 2006}, 
Panoramas et synth\`eses, SMF (2011), available on arXiv.org under 
{\tt http://arxiv.org/abs/0706.4447} 

\end{thebibliography}
\end{document}